\documentclass[10pt, letterpaper]{article}

\let\n\noindent

\newcommand{\beq}{\begin{equation}} \newcommand{\eeq}{\end{equation}}
\newcommand{\beqa}{\begin{eqnarray}}    \newcommand{\eeqa}{\end{eqnarray}}
\newcommand{\btab}{\begin{tabular}}     \newcommand{\etab}{\end{tabular}}

\newcommand{\bt}{\begin{table}}     \newcommand{\et}{\end{table}}

\newcommand{\ba}{\begin{array}}     \newcommand{\ea}{\end{array}}
\newcommand{\bc}{\begin{center}}        \newcommand{\ec}{\end{center}}
\newcommand{\bfig}{\begin{figure}}      \newcommand{\efig}{\end{figure}}
\newcommand{\bp}{\begin{picture}}       \newcommand{\ep}{\end{picture}}
\newcommand{\bq}{\begin{quote}}     \newcommand{\eq}{\end{quote}}
\newcommand{\ben}{\begin{enumerate}}    \newcommand{\een}{\end{enumerate}}

\font\tenmsy=msbm10
\font\sevenmsy=msbm10 at 7pt
\font\fivemsy=msbm10 at 5pt
\newfam\msyfam 
\textfont\msyfam=\tenmsy
\scriptfont\msyfam=\sevenmsy
\scriptscriptfont\msyfam=\fivemsy
\def\blackB{\fam\msyfam\tenmsy}
\def\Z{{\blackB Z}}

\def\frac#1#2{{\textstyle{#1\over #2}}}

\def\text#1{\quad\hbox{#1}\quad}

\def\la{\lambda}

\def\ka{\kappa}

\def\At{{\tilde{A}}}
\def\Jt{{\tilde{J}}}

\def\At{{\tilde{A}}}

\def\om{{\omega}}
\def\y{{\infty}}

\def\rw{\rightarrow}

\def\frac#1#2{{#1 \over #2}}

\def\rw{{\rightarrow}}

\begin{document}

\title{Jagged partitions}

\author{J.-F. Fortin, P. Jacob and P.
Mathieu\thanks{jffortin@phy.ulaval.ca, pjacob@phy.ulaval.ca,
pmathieu@phy.ulaval.ca} \\
\\ 
D\'epartement de physique, de g\'enie physique et d'optique,\\
Universit\'e Laval, \\
Qu\'ebec, Canada, G1K 7P4.
}

\date{October 2003}

\maketitle

\begin{abstract}

By jagged partitions we refer to an ordered collection of non-negative integers $(n_1,n_2,\cdots , n_m)$ with $n_m\geq
p$ for some positive integer $p$, further subject to some weakly decreasing conditions that prevent them for
being genuine partitions. The case analyzed in greater detail here corresponds to $p=1$  and the
following  conditions
$n_i\geq n_{i+1}-1$ and
$n_i\geq n_{i+2}$.  A number of properties for the corresponding partition function are derived, including rather
remarkable congruence relations. An interesting application of jagged partitions concerns the derivation of generating
functions for enumerating partitions with special restrictions, a point that is illustrated with various examples.
\end{abstract}

\newpage


\section{Introduction}

A new type of `partitions' (dubbed `jagged' for a reason to be explained shortly) recently arose in the context of a
conformal field-theoretical problem, namely the construction of a  quasi-particle basis for graded
 parafermionic theory with $\Z_{K}$ cyclic symmetry. 
These jagged partitions are ordered sequences of non-negative integers:
\beq \label{dede}
(n_1,n_2,\cdots , n_m)\qquad n_j\geq0,
\quad n_m\geq 1\, ,
\eeq
satisfying the weakly decreasing conditions
\beq\label{oror}
 n_j\geq n_{j+1}-1\qquad {\rm and}  \qquad  n_j\geq n_{j+2}\, .
\eeq
In their original appearance context, the jagged partitions were further subject to an exclusion condition:
$n_i \geq  n_{i+K-1} +1$ or $
n_i = n_{i+K-1} $ and $n_{i+1} =  n_{i+K-2}+2$,  with $K$ even. The generating function for these restricted jagged
partitions has been obtained in \cite{BFJM} while the corresponding one for  $K$ odd was found in
\cite{FJM}.
This  counting problem appeared to be a rather direct generalization of the enumeration of ordinary partitions
$(\la_1,\la_2,\cdots ,
\la_m)$ subject to the restriction
$\la_i\geq
\la_{i+k-1}+2$ solved by Andrews \cite{Andrr}.  

Our previous results strongly suggest that these new `partitions' (in their unrestricted versions) could have rather nice
properties. One of our goal in this work is to pinpoint some of them. In particular, using the previously obtained
generating function for $j(n)$, the number of jagged partitions of $n$, we derive a recurrence
relation together with a number of simple congruence properties for $j(n)$ itself, the most spectacular being that
$j(8n+7)\equiv 0$ mod 64. We also show that $j(n)$ is given by the Cauchy product of $p(m-n)$, the number of partitions
of
$n-m$, and $d(m)$, the number of partitions of $m$ into distinct parts, a result that entails an exact expression for
$j(n)$.  These results are presented in section 2. Moreover, by adapting to our case the Ramanujan's method for
obtaining the  generating function of $p(5n+4)$, we have derived a number of interesting generating functions for
$j(rn+s)$ with fixed $r$ and
$s$.  These results are reported in section 3.

Another aim of the present work is to explore the use of jagged partitions as a tool for enumerating standard partitions
satisfying special restrictions.  Let us explain in which sense that can be done.  There exists a simple bijection
between unrestricted jagged partitions (that is, vectors satisfying (\ref{dede}) and (\ref{oror})) and partitions subject
to a `difference-two condition at distance 2' :
$\la_i\geq
\la_{i+2}+2$ (the jagged partitions being simply augmented by the addition of a staircase) \cite{BFJM, FJM}. The
corresponding generating functions are consequently related in a simple way.  Here we stress that phrasing 
 the counting problem in terms  of jagged partitions induces an important simplification within our working framework
\cite{BFJM}. This method amounts to derive the generating function as the solution of a recurrence relation,  a technique
inspired from that of Andrews \cite{Andrr} for obtaining recurrence relations for restricted partitions. When formulated
in terms of jagged partitions, this method leads to a first-order recurrence relation, while its reformulation  directly
in terms of ordinary partitions leads to a third-order recurrence relation. This point is discussed in full detail in
section 4, where we also highlight a convenient pictorial tool for deriving recurrence relations.

This suggests that one
could similarly  count generalized jagged partitions to  obtain the generating
functions for  partitions with other interesting special restrictions.

But in order to see how jagged partitions could naturally be generalized, let us emphasis some salient features of the
definition (\ref{dede})-(\ref{oror}).
We see from (\ref{oror}) that the vector $(n_1,\cdots, n_m)$ is not a  genuine partition because 
the non-decreasing condition is not satisfied: an increment by one unit from $n_j$ to $n_{j+1}$ is allowed.
However a further increase by one unit from $n_{j+1} $ to $n_{j+2}$ is ruled out by the second condition in (\ref{oror}).
Note also that if the last entry is 1, zero entries are allowed. 
For instance, the set of jagged partitions of weight $\leq 7$ and length 5 (the length and the weight being respectively
the number of parts $m$ and the sum of the parts) is
\beq\label{exja}
\matrix{& (10101) &(20101)&(11101) &(30101) &(21101) &(11111)\cr
 & (12101)& (40101) 
&(31101) &(22101) &(21111) &(12111)\cr
&(21201) & (50101)& (41101) &(32101) 
&(31111) & (22111) \cr &(23101)& (21211)
&(12121) & (31201)& (22201)\,.\cr}
\eeq
Note that the
lowest-weight jagged partition of given length 
has the form
$(\cdots 01010101)$. The jagged 
character of this `ground-state partition' accounts for their name. 

In view of introducing hierarchies of novel jagged partitions,  it is appropriate to
introduce a more precise and, at the same time, more flexible terminology for those jagged partitions we have been discussing
so far. Given that they are characterized by their ground state, these can be conveniently called
$01$-partitions. The 01 notation indicates that these are the (pseudo) partitions built on the ground state whose period is
01. Natural extensions are thus $02$-partitions or $001$- partitions, etc.  In this terminology, ordinary partitions are 
1-partitions. In section 5, we present the generating functions for 02-, 012-, 001-partitions with fixed length and
apply them to special enumeration  problems for ordinary partitions.

\section{Basic properties of  01-partitions functions}

\n {\bf Definition 1.} A 01-partition of $n$ is a weakly decreasing sequence of non-negative integers $(n_1,n_2,\cdots)$
where $n=\sum_i n_i$, such that the last entry in non-zero and the conditions (\ref{oror}) are satisfied.

\n {\bf Definition 2.} The function $j(n)$ is the number of 01-partitions of $n$.

\n {\bf Theorem 3.} The generating function for the 01-partitions is 
\beq \label{genue}
J(q)= \sum_{n=0}^\y j(n)q^n= {(-q)_\y\over (q)_\y}\,.
\eeq

\n {\it Proof. } This follows from \cite{BFJM}, Eq. (5.17), with $z=1$, using 
$(a)_n:= (a;q)_n=(1-a)(1-aq)\cdots(1-aq^{n-1})$. Another proof is presented in section 4.

\n {\bf Corollary 4.} The partition function $j(n)$ satisfies the recurrence relation
\beq
j(n)=2\sum_{m\geq 1}(-1)^{m+1} j(n-{m^2})\,, \qquad n\geq 1\,.\eeq

\n {\it Proof. }  This follows  from the reexpression of (\ref{genue})  as (cf. \cite{AAR}, the equation following
(10.6.6))
\beq \label{genuer}
J(q)= {1\over \theta_4(q)}\, , 
\eeq
where
\beq\label{thet}
\theta_4(q)= 1+2\sum_{m\geq 1}(-1)^m q^{m^2}\, .\eeq 
One then has
\beq \label{unite}
1= J(q) {1\over J(q)}=  \sum_{n=0}^\y j(n)q^n\left (1+2\sum_{m\geq 1}(-1)^m q^{m^2}\right)\, , \eeq
from which the result follows.

\n {\bf Corollary 5.} The partition function $j(n)$ is related to $p(n)$, the number of partitions of $n$ and
$d(n)$ the number of  partitions of $n$ into distinct parts, by
\beq \label{jpdn}
j(n)= \sum_{m=0}^n p(n-m)d(m)\, .\eeq
\n {\it Proof. } The generating functions for $p(n)$ and $d(n)$ are respectively
\beq 
P(q)=\sum_{n=0}^\y p(n)q^n= {1\over (q)_\y}\; ,\qquad D(q)=\sum_{n=0}^\y d(n)q^n= {(-q)_\y}\, .
\eeq
We have thus
\beq
J(q)= D(q)P(q)\,,
\eeq
which implies ({\ref{jpdn}).

\n {\bf Remark 6.} Given that closed-form expressions are known for both $p(n)$ and $d(n)$,  due
to Hardy-Ramanujan-Rademacher (cf. theorem 5.1 in  \cite{Andr}) and Hau-Iseki-Hagis respectively (cf. Ex 5.3 therein),
this entails a closed-form expression for
$j(n)$.  We recall that at the origin of the Hardy-Ramanujan collaboration on partition problems is rooted in the
fruitful `Ramanujan's false statement' (\cite{Har} p. 9) that the coefficient of $q^n$ in $1/\theta_4(q)$, which is
precisely
$j(n)$, is given by
\beq{1\over 4n}\left (\cosh \pi {\sqrt n} - {\sinh \pi {\sqrt n}\over \pi {\sqrt n}}\right)\, .\eeq
 This remains a remarkable approximation.
  
The expression of the generating function $J(q)$ in terms of the inverse of $\theta_4$ leads immediately to two simple
congruence properties.

\n {\bf Corollary 7.} $j(n)\equiv 0 $ mod 2 for $n>1$.

\n {\it Proof. } Eq. (\ref{unite}) can be written as
\beq \label{unitee}
1= J(q)+2 J(q)\sum_{n\geq 1}(-1)^{n} q^{n^2}\,, \eeq
which implies (using (\ref{genue}))
\beq\label{iter} 
J(q)= 1+2J(q)\sum_{n\geq 1}(-1)^{n+1} q^{n^2}=  1+2\, {(-q)_\y\over (q)_\y} \sum_{n\geq 1}(-1)^{n+1} q^{n^2}\, .
\eeq
01-partition functions are thus always even except for $j(0)=1$.

\n {\bf Corollary 8.} $j(n)\equiv 0 $ mod $2^{p+1}$ for all $n$ that cannot be written as a sum of $p$ or less squares.

\n {\it Proof. } Iterating  the first equality in (\ref{iter}), we obtain
\beq\label{itera} 
J(q)= 1+\sum_{p=1}^\y 2^p\sum_{n_1,n_2,..n_p\geq 1}(-1)^{\sum_{i=1}^p (n_i+1)}\,  q^{\sum_{i=1}^p n_i^2}\,,
\eeq
from which the statement is immediate.


This result is meaningful only for $p\leq 3$ since, by  Lagrange theorem, all numbers can be written as a sum of four
squares. Here is  a simple application of the corollary: since $(8n+7)$ cannot be written as a sum of less that 4
squares,  it readily follows that $j(8n+7) \equiv 0 $ mod
16.

\n {\bf Corollary 9.} If $p'$ stands for the least number of squares into which $rn+s$, with $1\leq s<r$, is
decomposable, then
$j(rn+s)\equiv 0$ mod $a\, 2^{p'}$, where $a= {\rm min}\, (c, 2)$ with $c$ being the number of distinct vectors
$(n'_1,\cdots, n'_{p'})$ of strictly positive entries which sum to $s$.

\n {\it Proof. } The least number of squares into which a number can be decomposed is easily found by inspection of the
non-zero values of $m^2$ modulo $r$, which are 
$ (1,4,\cdots, (r-1)^2)$ mod $r$. The number of solutions of $s=n'_1+\cdots+ n'_{p'}$, where the $n'_i$ are square
residues mod $r$ (i.e., $1\leq n'_i< r$), is a divisor of  the number of equivalent terms in the first summation in
(\ref{itera}). For $c>2$, the $p'$-square contribution might not be the lowest one and  we need to consider also the
contribution of $p'+1$ squares; therefore the bound obtained by the lowest contributing term is actually
$a= {\rm min}\, (c, 2)$.

Let us consider some applications of this corollary.  An immediate consequence is that $j(rn+2)\equiv 0$ mod 4 when $r>2$.
Indeed, these numbers can be decomposed in two squares, not less, hence $p'=2$, and this can be done in one way mod $r$
$(2=1+1)$, so that
$c=1$. Similarly, $j(rn+3)\equiv 0$ mod 4 when $r>3$ since then $p'=3$ and the combinatorial factor $c$ is again 1. In
the same vein, it follows that $j(rn+5)\equiv 0$ mod 8 when $r>5$ since $rn+5$ can be decomposed as two squares and
$s=4+1=1+4$ so that $c=2$.  Our final example deserves to be singled out.

\n {\bf Corollary  10.} $j(8n+7)\equiv0 \; {\rm mod}\; 64$.

\n {\it Proof. }  We
have already noticed that the first coefficient of (\ref{itera}) contributing to $j(8n+7)$ comes from the
decomposition into 4 squares. But it has a combinatorial factor
 at least equal to 4 since $7=4+1+1+1$ and its three
permutations. It is larger than 4 whenever the three residues $n_i'=1$ are residues of distinct numbers $n_i$. On the
other hand,  since
$j(8n+7)$ cannot be decomposed into five squares, the factor
$a$ in Corollary 9 can be replaced by
$a'= {\rm min}\, (c, 4)= 4$.

Here is an illustration of this congruence from a direct computation based on (\ref{itera}):
\beq
j(15)= -2^4\cdot 12-2^6\cdot 20+2^7\cdot 7+2^9\cdot 36-2^{12}\cdot 12+2^{15}=64\cdot 23\, .\eeq
The number of contributing terms from the first summation  is 12 since there are 12 ways of writing
$15=\sum_{i=1}^4n_i^2= 9+4+1+1$ and permutations; modulo 8, this equation reduces to $7=4+1+1+1$ and permutations;
hence here there are three $n'_i$ equal to 1  but the corresponding $n_i$'s are all different (they are permutations
of
$(9,4,1)$).

\section{Ramanujan-type generating functions for  01-partitions}

In this section, we derive special generating functions analogous to the famous Ramanujan result (cf. \cite{Andr} chap
10):
\beq \label{rama}
\sum_{n=0}^\y p(5n+4)q^n= 5{(q^5;q^5)^5_\y\over (q)_\y^6}\,,\eeq
which implies  that $p(5n+4)\equiv 0$
mod 5.

We have calculated  the generating functions for all 01--partitions $j(rn+s)$ with $2\leq r\leq 8$.
Those for $r=2,3,4$ take a rather compact form:

\n {\bf Proposition 11.} The generating functions for the 01-partitions $j(rn+s)$, $2\leq r\leq 4 $ read 
\begin{eqnarray}
&&  \sum_{n=0}^\y j(2n)q^n = {(q^4;q^4)_\y^5\over (q)_\y^4(q^8;q^8)_\y^2}\; , \qquad\qquad
\sum_{n=0}^\y j(2n+1)q^n = 2 \,{(q^2;q^2)_\y^2 (q^8;q^8)_\y^2\over (q)_\y^4(q^4;q^4)_\y}\;,\cr
&&
\sum_{n=0}^\y j(3n)q^n = {(q^2;q^2)_\y^4(q^3;q^3)_\y^6\over (q)_\y^8(q^6;q^6)_\y^3}\; , \qquad\qquad
\sum_{n=0}^\y j(3n+1)q^n = 2\,{(q^2;q^2)_\y^3(q^3;q^3)_\y^3\over (q)_\y^7}\;,\cr
&&\sum_{n=0}^\y j(3n+2)q^n = 4\,{(q^2;q^2)_\y^2(q^6;q^6)_\y^3\over (q)_\y^6}\; , \qquad\qquad
\sum_{n=0}^\y j(4n)q^n ={(q^2;q^2)_\y^{19}\over (q)_\y^{14}(q^4;q^4)_\y^6}\;,\cr 
&&\sum_{n=0}^\y j(4n+1)q^n = 2\,{(q^2;q^2)_\y^{13}\over (q)_\y^{12}(q^4;q^4)_\y^2}\; , \qquad\qquad
\sum_{n=0}^\y j(4n+2)q^n = 4\,{(q^2;q^2)_\y^7(q^4;q^4)_\y^2\over (q)_\y^{10}}\;,\cr
&&\sum_{n=0}^\y j(4n+3)q^n = 8\,{(q^2;q^2)_\y(q^4;q^4)_\y^6\over (q)_\y^8}\;.
\end{eqnarray} 

Note that from these expressions, we obtain three more formulas for $J(q)$, the simplest one being: 
\beq
J(q)= \sum_{n=0}^\y [j(2n)q^{2n}+ j(2n+1)q^{2n+1}]= {(q^8;q^8)_\y^5\over (q^2;q^2)_\y^4(q^{16};q^{16})_\y^2}+
2 q\,{(q^4;q^4)_\y^2 (q^{16};q^{16})_\y^2\over (q^2;q^2)_\y^4(q^8;q^8)_\y}\;.
\eeq

The last result we make explicit is the following, which provides an alternative proof of Corollary 10.

\n {\bf Proposition 12.} The generating functions for the 01-partitions $j(8n+7)$ is
\beq
\sum_{n=0}^\y j(8n+7)q^n =
64\,{(q^2;q^2)_\y^{22}\over (q)_\y^{23}}\;.
\eeq


\n {\it Sketch of the proof of Propositions 11 and 12.} Our proof of these relations is inspired by the derivation of
(\ref{rama}) given in
\cite{AAR} Sect. 11.7. The proof is divided into three steps.

\n {\bf 1}: We first compute the ratio $J(q^{1/r})/J(q^r)$. We keep the denominator in the form of a product and break 
the inverse of $J(q^{1/r})$,
\beq \label{deder}
{1\over J(q^{1/r})}= \sum_{n=-\y}^\y (-1)^n q^{n^2/r}\, ,\eeq
 into $r$ sums (setting $n=rm+\ell$, with $0\leq \ell\leq r-1$) and evaluate each sum by means of the Jacobi
triple-product identity:
\beq \sum_{n=-\y}^\y (-1)^n q^{n(n+1)/2}z^n= (qz)_\y\, (z^{-1})_\y\,(q)_\y\,.
\eeq
 That leads us to an
expression of the form:
\beq\label{debut}
{J(q^{1/r})\over J(q^r) } = {1\over \sum_{\ell=0}^{r-1}\xi_\ell(q) q^{\ell/r}}\,.
\eeq
The different $\xi_\ell(q)$'s are now known functions of $q$.  The division by
$J(q^r)$ is purely conventional, its role being to simplify the form of the `coefficients' $\xi_\ell(q)$. In
particular, for $r$ odd, it makes
$\xi_0(q)=1$.  By construction, the $\xi_\ell(q)$'s are invariant under the transformation $q^{1/r}\rw q^{1/r}\om^k$, with
$\om= e^{2\pi i/r}$.

\n {\bf 2}: We then 
observe that by rewriting 
\beq \label{nuden}
{1\over \sum_{\ell=0}^{r-1}\xi_\ell(q)\,  q^{\ell/r}} ={\prod_{k=1}^{r-1}\sum_{\ell=0}^{r-1}\xi_\ell(q)
\, \om^{k\ell}\, q^{\ell/r}\over
\prod_{k=0}^{r-1}\sum_{\ell=0}^{r-1}\xi_\ell(q)\,  \om^{k\ell}\, q^{\ell/r}}=:
{\prod_{k=1}^{r-1}\sum_{\ell=0}^{r-1}\xi_\ell(q)\, 
\om^{k\ell}\, q^{\ell/r}\over \Omega_r(q)}\,,
\eeq 
the denominator becomes a function of $q$ and not of $q^{1/r}$. This denominator is evaluated from (\ref{debut}), but 
using this times the product form of both $J(q^{1/r})$ and $J(q^{r})$:
\beq \label{defiomeg}
{1\over \Omega_r(q)} = 
{1\over [J(q^r) ]^r} \prod_{k=0}^{r-1} J(\om ^k q^{1/r})\,.
\eeq
The evaluation of $\Omega_r(q)$ relies on the lemma that follows this sketch.

\n {\bf 3}: 
Finally, in order to select the terms $j(rn-s')$ from $J(q^{1/r})$, we multiply it by $q^{s'/r}$, replace  $q^{1/r}$
by
$\om^kq^{1/r}$ and then sum over all values of $k$.  For $J(q^{1/r})$, we use formula (\ref{debut}). The same operation of
the rhs selects the $q^{s'/r+m}$ terms in the numerator (with $m$ integer). The final result follows by multiplying both
sides by
$J(q^r)$.

\n {\bf Lemma 13.}  The  product $\prod_{k=0}^{r-1} J(\om ^k q^{1/r})$, for $r$ odd and  prime or $r=2^p$ takes the
following form:
\begin{eqnarray}\label{lelem}
\prod_{k=0}^{r-1}\prod_{n=1}^{\y} {1+(\om^kq^{1/r})^n\over 1-(\om^kq^{1/r})^n}&=& {(q^2;q^2)_\y^{r+1} \,
(q^r;q^r)^2_\y\over (q)_\y^{2r+2}\, (q^{2r};q^{2r})_\y }\qquad  {\rm for}æ\quad r\; {\rm odd ~and~ prime}\;,\cr
&=& {(q^{2};q^{2})_\y^{r}\over (q)_\y^{2r} } \qquad {\rm for}æ\quad r=2^p \;.
\end{eqnarray}

\n {\it Proof. }  Set  $n=rm+s$ with $0<s\leq r$ and $m\geq 0$.  In the case $r=2^p$, for given values of $m$ and 
$s\not=r$, the
$r$ terms (resulting from the product over $k$) of the numerator are easily seen to be permutations of those of the
denominator so that the product is 1; for 
$s=r$, the $k$ dependence disappear (all $r$ terms are identical) and the product over $m$ yields
$(-q)_\y^r/(q)_\y^r$. For
$r$ prime and $\not=2$, we find that for $m$ fixed:
\beq
\prod_{s=1}^{r}\prod_{k=0}^{r-1}{1+(\om^kq^{1/r})^{rm+s}\over 1-(\om^kq^{1/r})^{rm+s}}=
{(1+q^{rm+1})\cdots (1+q^{rm+r-1}) (1+q^{m+1})^r 
\over (1-q^{rm+1})\cdots (1-q^{rm+r-1}) (1-q^{m+1})^r}\cdot{(1+q^{rm+r})(1-q^{rm+r})\over
(1-q^{rm+r})(1+q^{rm+r})}\;,
\eeq
(where at the end we have introduced a suitable decomposition of 1). The quoted result is obtained by taking the product
over all
$m$ and reorganizing the sums using repeatedly the simple identity:
\beq
(-q^c;q^c)_\y= {(q^{2c};q^{2c})_\y\over (q^c;q^c)_\y}\,.\eeq
These cases cover all those (namely $r=3,\, 2^1,\,2^2$ and $2^3$) needed to tackle Propositions 11 and 12.

\n {\it Proof of Proposition 11 for $r=3$.} Let us now detail the case $r=3$. By decomposing the sum in (\ref{deder})
into three sums according to $n=3m, 3m\pm1$ and transform them by the Jacobi triple-product identity:
\begin{eqnarray}
 {1\over J(q^{1/3})}&=&  \sum_{m=-\y}^\y (-1)^m q^{3m^2}-2q^{1/3}\sum_{m=-\y}^\y (-1)^m q^{3m^2+2m}\cr
&=& {(q^3;q^3)_\y^2\over (q^6;q^6)_\y} -2q^{1/3} {(q)_\y \, (q^6;q^6)_\y^2\over  (q^2;q^2)_\y\, (q^3;q^3)_\y}\,,
\end{eqnarray}
we end up with
\beq \label{potu}{ J(q^{1/3})\over J(q^3)}= {1\over 1-2q^{1/3}\zeta_1(q) }\,,
\eeq
where
\beq \zeta_1(q)= {(q)_\y \, (q^6;q^6)_\y^3\over  (q^2;q^2)_\y\, (q^3;q^3)^3_\y}\,.
\eeq
The function $\zeta_1(q)$   is manifestly unaffected by the replacement  of  $q^{1/3}$ by $ q^{1/3}\om^k$, with
$\om= e^{2\pi i/3}$. We thus have $\xi_0(q)=1\,, \, \xi_2(q)=0$ and we have  redefined $\xi_1(q)=-2\zeta_1(q)$.
This completes step 1.  We now evaluate $\Omega_3(q)$:
\beq \Omega_3(q)= \prod_{k=0}^2 [{1-2q^{1/3}\om^k \zeta_1(q) }] = 1-8q\zeta_1^3(q)\,,\eeq
while, from (\ref{defiomeg}) and (\ref{lelem}), we have
\beq {1\over \Omega_3(q)}= {(q^2;q^2)^4_\y, (q^3;q^3)_\y^8\over (q)_\y^8 (q^6;q^6)^4_\y}\,.\eeq
We can thus write
\beq
{1\over 1-2q^{1/3}\zeta_1(q) } = {\prod_{k=1}^{2}\left[1-2q^{1/3}\om^k \zeta_1(q)\right] \over
1-8q\zeta_1^3(q)}= {1+2q^{1/3}\zeta_1(q)+4q^{2/3}\zeta_1^2(q)\over   \Omega_3(q)}\,,\eeq
and by (\ref{potu}), this is equal to $J(q^{1/3})/J(q^3)$:
\beq \label{uty}
J(q^{1/3})= \sum_{n=0}^\y j(n)\, q^{n/3}=  {J(q^3)\over \Omega_3(q)}
\left(1+2q^{1/3}\zeta_1(q)+4q^{2/3}\zeta_1^2(q)\right)\,.\eeq
Let us now replace  $q^{1/3}$ by $ q^{1/3}\om^k$ on both sides and sum the resulting equation over $k$. On the lhs, we get
\beq
\sum_{k=0}^2 \sum_{n=0}^\y j(n)\, (q^{1/3}\om^k)^n=  \sum_{n=0}^\y j(n)\,q^{n/3}(1+\om^n+\om^{2n})=3 \sum_{m=0}^\y
j(3m)\,q^{m}\,,\eeq
while on the rhs, we find
\beq
\sum_{k=0}^2 \left(1+2q^{1/3}\, \om^k\, \zeta_1(q)+4q^{2/3}\,\om^{2k}\, \zeta_1^2(q)\right)= 3\,,\eeq
so that
\beq
 \sum_{m=0}^\y
j(3m)\,q^{m}= {(q^2;q^2)_\y^4\, (q^3;q^3)_\y^6\over (q)_\y^8\, (q^6;q^6)_\y^3}\,.\eeq
By multiplying (\ref{uty}) by $q^{1/3}$, replacing  $q^{1/3}$ by $q^{1/3}\om^k$ on both sides and summing over $k$, we get
\beq
 \sum_{m=0}^\y
j(3m-1)\,q^{m}= q\sum_{m=0}^\y
j(3m+2)\,q^{m}= 4q\, {(q^2;q^2)_\y^2\, (q^6;q^6)_\y^3\over (q)_\y^6}\,.\eeq
The generating function for $j(3n+1)$ is obtained in a similar way, by multiplying (\ref{uty}) by $q^{2/3}$,  
$q^{1/3}\rw q^{1/3}\om^k$ and sum over $k$.

\n {\it Proof of Proposition 12.} It is convenient to start with $J(q^{1/4})$:
\begin{equation}
\frac{1}{J(q^{1/4})}=\sum_{n=-{\infty}}^{\infty}{(-1)^nq^{n^2/4}}
=\sum_{n=-{\infty}}^{\infty}q^{4n^2}-2q^{1/4}
\sum_{n=-{\infty}}^{\infty}q^{2n(2n+1)}+q\sum_{n=-{\infty}}^{\infty}q^{4n(n+1)}\,.
\end{equation}
Using the Jacobi triple-product identity, this becomes
\begin{equation}
\frac{1}{J(q^{1/4})}=\frac{(q^8;q^8)_{\infty}^5}
{(q^4;q^4)_{\infty}^2(q^{16};q^{16})_{\infty}^2}
-2q^{1/4}\frac{(q^4;q^4)_{\infty}^2}
{(q^2;q^2)_{\infty}}+2q\frac{(q^{16};q^{16})_{\infty}^2}
{(q^8;q^8)_{\infty}}\,.
\end{equation}
Therefore ${J(q^8)}/{J(q^{1/8})}$ can be written as:
\begin{equation}
\frac{J(q^8)}{J(q^{1/8})}=\zeta_0-2q^{1/8}\zeta_1+2q^{1/2}\zeta_4\,,
\end{equation}
with
\begin{equation} \label{leszeta}
\zeta_0=\frac{(q^4;q^4)_{\infty}^5(q^{16};q^{16})_{\infty}}
{(q^2;q^2)_{\infty}^2(q^8;q^8)_{\infty}^4}\;,
\qquad 
\zeta_1=\frac{(q^2;q^2)_{\infty}^2(q^{16};q^{16})_{\infty}}
{(q;q)_{\infty}(q^8;q^8)_{\infty}^2}\;,
\qquad
\zeta_4=\frac{(q^{16};q^{16})_{\infty}} {(q^4;q^4)_{\infty}}\,.
\end{equation}
Using (\ref{nuden}), we can write
\begin{equation}
\frac{J(q^{1/8})}{J(q^{8})}=
\frac{q^{7/8}\alpha(\zeta_0,\zeta_1,\zeta_4)+\cdots }{\Omega_8(q)}\,,
\end{equation}
where we only write the coefficient of the terms $q^{7/8+m}$ for $m$ integer, which is
\begin{equation}
\alpha(\zeta_0,\zeta_1,\zeta_4)=64
(2\zeta_1^7-\zeta_0^3\zeta_1^3\zeta_4-4q\zeta_0\zeta_1^3\zeta_4^3)\,,
\end{equation}
while
\begin{equation}
\frac{1}{\Omega_8(q)}=\frac{(q^2;q^2)_{\infty}^8(q^8;q^8)_{\infty}^{16}}
{(q;q)_{\infty}^{16}(q^{16};q^{16})_{\infty}^8}\,.
\end{equation}
Note that ${\Omega_8(q)}$ is actually a function of $q$.
For $j(8n+7)$, we have thus:
\begin{equation}
\frac{(q^8;q^8)_{\infty}^2}{(q^{16};q^{16})_{\infty}}\sum_{n=0}^{\infty}j(8n+7)q^n=
64(2\zeta_1^7-\zeta_0^3\zeta_1^3\zeta_4-4q\zeta_0\zeta_1^3\zeta_4^3)
\frac{(q^2;q^2)_{\infty}^8(q^8;q^8)_{\infty}^{16}}
{(q;q)_{\infty}^{16}(q^{16};q^{16})_{\infty}^8}\,.
\end{equation}
The result follows from the identity
\beq \label{idrty}
\zeta_0^3\zeta_1^3\zeta_4+4q\zeta_0\zeta_1^3\zeta_4^3 = \zeta_1^7\,,
\eeq
which is established in the appendix.

\section{Generating function for  01-partitions of prescribed length}


In this section, we first rederive the generating function for 01-partitions by using a modified version of a
technique  introduced in App. A of \cite{BFJM} .  We first motivate the method by deriving a recurrence relation for
ordinary partitions
$p(m,n)$ of $n$ into $m$ parts in the language of Ferrer graphs. Every partition of length $m$ of
$n$  can be represented by graphs of $n$ dots distributed among $m$ rows whose length (specified by the parts of the
partition) does not increase when read from top to bottom. The set of such Ferrer graphs of $m$ rows can be decomposed into
the disjoint union of those graphs with precisely one dot on the last row and those with more that one dot
on this
$m$-th row. The former set is characterized by the fact that removing the single dot on the $m$-th row leaves a
Ferrer graph with $m-1$ rows (and of course a total of $n-1$ dots). Similarly, the graphs with two dots or more on
the $m$-th row are characterized by the fact that they still have $m$ rows if we take out the first column (reducing
$n$ to
$n-m$). That leads directly to the following recurrence relation
\beq \label{parsisi}
p(m,n)=p(m-1,n-1)+p(m,n-m)\,.
\eeq

Instead of working with graphs, we will introduce a pictorial representation that captures the notion of `building on a
ground state with prescribed structure'. Let us then represent
$p(m,n)$ as 
\beq
p(m,n): \, (\cdots 11111^+) \qquad (m\;{\rm entries)}\,,
\eeq
where the symbol $+$ indicates that we can build up the partition on the indicated ground state up to the position of  the
$+$.
It is clear that this `filling process' can be broken as follows:
\beq \label{pictu}
(\cdots 11111^+)= (\cdots 1111^+1)+ (\cdots 22222^+)\,,
\eeq
(from now on, all vectors are supposed to have $m$ entries). In the first term of the left hand side, we isolate 
those partitions that have at least one 1. What is left is completely taken into account by the
set $(\cdots 22222^+)$, which describes all partitions without any 1.  But this pictorial relation is nothing but the
Ferrer graph retranscription of the decomposition considered above.   Indeed, elements of
the set
$(\cdots 22222^+)$ are in a one-to-one correspondence with those enumerated by $p(m,n-m)$: the correspondence is obtained 
by subtracting the partition $(1^m)$ from elements of the former set.  Similarly, there is a one-to-one correspondence
between elements of
$(\cdots 1111^+1)$ and the partitions of length $m-1$ and weight $n-1$: one simply adds a 1 at the extremity of the
elements of the latter set, an operation that preserves the non-increasing  character of the partition. Therefore
(\ref{pictu}) is equivalent to  (\ref{parsisi}).

The recurrence  relation (\ref{parsisi}) is far from being new: it has first been obtained by Euler
(see e.g., \cite{Rad} Eq. (97.7) but quoted with a misprint) from the generating function $P_m(q)= \sum_n p(m,n)q^n=
q^m/(q)_m$. The relation follows from the identity $(1-q^m)P_m= qP_{m-1}$. Euler used it to compute
recursively the
$p(m,n)$'s. Our aim is just the opposite: we want to find generating
functions by solving recurrence relations.

Our `pictorial' derivation of (\ref{parsisi}) is inspired by the method used by Andrews to find the
generating function for partitions with prescribed length and `difference-two at distance $k-1$' \cite{Andrr}. The idea is
to split the set we want to enumerate into distinct sets characterized by their boundary condition. The method is here
adapted to the case where the restriction is washed out.

\n {\bf Definition 14.} The function $j(m,n)$ yields the number of 01-partitions of $n$ having exactly $m$ parts (i.e.,
of length
$m$).

\n {\bf Lemma 15.} The function $j(m,n)$ satisfies the recurrence relations:
\begin{eqnarray} \label{petite}
j(m,n)& = &  j(m-2,n-1)+k(m,n) \cr
k(m,n)& = &  k(m-1,n-1)+j(m,n-m) \,,
\end{eqnarray}
where $k(m,n)$ stands for the number of 01-partitions of $n$ having exactly $m$ parts $\geq 1$ (i.e., having  no 0).

\n {\it Proof.} The functions  $j(m,n)$ and  $k(m,n)$
have the following pictorial representations:
\beq
j(m,n): \, (\cdots 010101^+) \; , \qquad\qquad  k(m,n): \, (\cdots 111111^+) \,.
\eeq
(This last  set should not be confused with the set of ordinary partitions of $n$ of length $m$; $k(m,n)$ really counts
special 01-partitions: for instance (1212) is an element of the set whose cardinality is $k(4,6)$.)
 The following relations are satisfied
\begin{eqnarray} 
(\cdots 010101^+) & = & (\cdots 0101^+01)+ (\cdots 111111^+)\cr
(\cdots 111111^+) & = & (\cdots 11111^+1)+ (\cdots 121212^+)\,.
\end{eqnarray} 
In the first relation, we isolate the terms with at least one pair 01 at the end, i.e., containing at least one 0. What
is left is the set of 01-partitions without parts equal to 0, hence built on $(\cdots 111)$. The second relation is
similar. These are equivalent to the recurrence relations (\ref{petite}). Let us explain the last term: the jagged
partitions
$(\cdots 121212^+)$ can be characterized by the following property: subtracting from these elements the partition
$(1^m)$ (reducing the weight by
$m$ but without affecting the length)  makes them elements of the set enumerated by
$j(m,n-m)$.  The functional retranscription of the other terms is obvious.

\n {\bf Theorem 16. } The generating function of $j(m,n)$ is 
\beq \label{genfu}
J(z;q)= {(-zq)_\y\over (z^2q)_\y} \,.
\eeq

\n {\it Proof. } 
Let us introduce the two generating functions
\beq
J(z;q)= \sum_{m,n=0}^\y j(m,n)z^mq^n\;, \qquad K(z;q)= \sum_{m,n=0}^\y k(m,n)z^mq^n\,.
\eeq
The recurrence relations (\ref{petite}) are then lifted to
\begin{eqnarray} 
J(z;q)& = &  z^2qJ(z;q)+K(z;q) \cr
K(z;q)& = &  zqK(z;q)+J(zq;q) \,.
\end{eqnarray}
Since the $q$-dependence is never affected, it is usually omitted. The solution of the above relations is easily obtained
\beq 
J(z)=  {K(z)\over 1-z^2q}= {J(zq)\over  (1-zq)(1-z^2q)}=\cdots =  {(-zq)_\y\over (z^2q)_\y}\;, \quad \qquad 
 K(z)=  {(-zq)_\y\over
(z^2q^2)_\y}\,,
\eeq
(where the dots stand for the infinite  iteration  of the preceding result).

Note that in contrast with $p(m,n)$, the function $j(m,n)$ does not necessarily vanish if $n<m$. 
For instance $j(4,2)=1$
and the corresponding 01-partition is $(0101)$.  As a further illustration, the $q^3$ coefficient of $J(z)$
obtained by expanding (\ref{genfu}) and the corresponding 01-partitions are:
\beq
z+2z^2+2z^3+z^4+z^5+z^6 : \quad\{(3)\,,(21)\,,(12)\,,(201)\,,(111)\,,(1101)\,,(10101)\,,(010101)\}\,.\eeq

The following lemma presents an application of the previous result to the enumeration of 01-partitions with at most $m$
parts.  Recall that  $p_m(n)$, the number
of partitions of $n$ into at most $m$ parts, is equal to $p(m,n+m)$. The analogous relation in the jagged
case is given in the following lemma.

\n {\bf Lemma 17. } The partition function $j_m(n)$ that counts the number of  01-partitions of $n$  with at most
$m$ parts satisfies
\beq\label{petit}
j_m(n)= j(m,n+m) -j(m-2,n+m-1)= k(m,n+m)\,,
\eeq
where $k(m,n)$ is defined in Lemma 15.

\n {\it Proof. }
We  have 
\beq
j(m,n+m)= j_m(n)+j(m-2,n+ m-1)\,.
\eeq
Indeed, by adding the partition $(1^m)$ to those counted by $j_m(n)$ we obtain all 01-partitions of weight
$n+m$ having no 0. To generate the whole set of 01-partitions of weight $n+m$ we simply need to add to this set  all those
01-partitions that contain at least one 0, that is, one pair of 01. By stripping this tail 01, we see that all these elements
are in correspondence with those counted by $j(m-2,n+m-1)$. This yields the above relation. 
The derivation makes clear that $j_m(n)$ is exactly the set $k(m,n+m)$. Moreover, (\ref{petit}) is simply the first
relation in (\ref{petite}).

Let us recall the bijection introduced in \cite{FJM}: by adding the staircase $(m-1,m-2,\cdots,1,0)$ to the vector
$(n_1,\cdots n_m)$, we transform it into an ordinary partition. With $\la_j=n_j+m-j$, the weakly decreasing conditions
(\ref{oror}) become
\beq \label{parc}
\la_j\geq \la_{j+1} \qquad {\rm and} \qquad \la_j\geq \la_{j+2}+2\,.
\eeq
To transform a generating function for  jagged partitions to one for partitions subject to (\ref{parc}), we 
 replace $z^N$ by $z^N q^{N(N-1)/2}$ within the sum defining the generating function.  We thus reproduce the result
of Andrews \cite{Andrr} (his $F_{3,3}(z;q)$) as follows (cf. Corollary 9 of
\cite{FJM}):

\n {\bf Corollary 18. } The generating function for the number of partitions satisfying
$\la_j\geq \la_{j+2}+2$ is  given by 
\beq
\sum_{m,n\geq 0} j(m,n) z^mq^{m(m-1)/2+n}= \sum_{m_0,m_1\geq 0}
{q^{(m_0+m_1)^2+m_1^2}\, z^{m_0+2m_1}\over (q)_{m_0}(q)_{m_1} }\;.
\eeq

\n {\bf Remark  19. }  Let us argue that the jagged-partition formulation of the problem is efficient by reworking the
counting  problem directly at the level of partitions.  Let us then introduce three sets of partitions  subject to
(\ref{parc}):
\beq
a(m,n): \, (\cdots 553311^+) \; ,\qquad  
b(m,n):   \, (\cdots 55331^+)\;,\qquad  
c(m,n):   \,
(\cdots 765432^+)\,.
\eeq
These sets obey the recurrence relations
\begin{eqnarray} 
a(m,n)&=& a(m,n-m)+b(m-1,n-1)\cr
b(m,n)&=& a(m-1,n-2m+1)+c(m,n)\cr
c(m,n)&=& a(m,n-2m)+c(m-1,n-m-1)\,.
\end{eqnarray}
For their generating functions (defined in the obvious way), this becomes
\begin{eqnarray} 
A(z)&=& A(zq)+zqB(z)\cr
B(z)&=& zqA(zq^2)+C(z)\cr
C(z)&=& A(zq^2)+zq^2C(zq)\,.
\end{eqnarray}
Solving for $A(z)$, we end up with
\beq
A(z)= (1+zq)A(zq)+z^2q^2A(zq^2)-z^3q^5A(zq^3)\,.\eeq
Although there are ways to solve this relation (using suitable transformations), ending up with a third-order
$q$-equation instead of a first-order one is enough to illustrate our point:  reformulating the problem in
terms of jagged partitions bear some magic simplifications.


\section{Generating function for  generalized jagged partitions of prescribed length}

We now consider three types of generalized jagged partitions, for which we also find  generating functions by solving a
first-order recurrence relation.  Since the construction of these generating functions  follow the same pattern as
for the 01-case, we condense the presentation by including the derivations of the recurrence relations within the proof
of the theorems in the first two cases and limit ourself to  displaying the
key recurrence relation and the final result in the last case.

\subsection{02-partitions}

\n {\bf Definition 20.} $02$-partitions  are defined as vectors $(n_1,\cdots, n_m) $ of
non-negative entries
$n_j\geq 0$,with
$n_m\geq 2$, subject to the  weak ordering conditions:
\beq 
n_j\geq n_{j+1}-2\;,\qquad  n_j\geq n_{j+2}\,.
\eeq



\n {\bf Theorem 21.} The generating function for ${\cal J}(m,n)$, the number of $02$-partitions of $n$ with length $m$,
is 
\beq \label{casdeu}
 {\tilde {\cal J}}(z)=\sum_{n,m\geq0} {\cal J}(m,n)z^mq^n= {(z^3q^6;q^3)_\y \over (zq^2)_\y
(z^2q^2)_\y}\,.
\eeq

\n {\it Proof. } 
The first step of the proof amounts to derive appropriate recurrence relations.
We introduce the following pictorial representation (function-set bijection) for
the 
${\cal J}(m,n)$ and two auxiliary functions ${\cal K}(m,n)$ and ${\cal L}(m,n)$
\beq 
{\cal J}(m,n): \,(\cdots 020202^+)\; , \qquad {\cal K}(m,n): \, (\cdots 121212^+) \; , \qquad
 {\cal L}(m,n):\, (\cdots 222222^+)\,.
\eeq
again with the understanding that all vectors have $m$ components.
While ${\cal J}(m,n)$
counts all 
$02$-partitions, ${\cal K}(m,n)$ counts those with no 0
and 
${\cal L}(m,n)$ enumerates those having neither 0 nor 1. The
following  relations hold between the different  `ground-state fillings':
\begin{eqnarray}
& & (\cdots 0202^+)=(\cdots 02^+02) + (\cdots 1212^+)\cr
& & (\cdots 1212^+)=(\cdots 12^+12) + (\cdots 2222^+)+ (\cdots 13^+13) \cr
& & (\cdots 2222^+)=(\cdots 222^+2) + (\cdots 2323^+)\,.\ 
\end{eqnarray}
These translate into the recurrence relations:
\begin{eqnarray}
& & {\cal J}(m,n)= {\cal J}(m-2,n-2) +{\cal K}(m,n)\cr
& & {\cal K}(m,n)= {\cal K}(m-2,n-3) +{\cal L}(m,n)+{\cal J}(m-2,n-m-2)\cr
& & {\cal L}(m,n)= {\cal L}(m-1,n-2) +{\cal K}(m,n-m)\,.\
\end{eqnarray}
As usual, these recurrence relations are then lifted to $q$-difference equations for their generating
functions
and these  take the form
\begin{eqnarray}
& & {\tilde {\cal J}}(z)= z^2q^2{\tilde {\cal J}}(z) +{\tilde {\cal K}}(z)\cr
& & {\tilde {\cal K}}(z)= z^2q^3{\tilde {\cal K}}(z) +{\tilde {\cal L}}(z)+z^2q^4{\tilde {\cal J}}(zq)\cr
& & {\tilde {\cal L}}(z)= zq^2{\tilde {\cal L}}(z) +{\tilde {\cal K}}(zq)\,.\
\end{eqnarray}
The auxiliary functions ${\tilde {\cal K}}(z)$ and ${\tilde {\cal L}}(z)$ can be eliminated and we get
\beq \label{redt}
{\tilde {\cal J}}(z)= {(1+zq^2+z^2q^4)\over (1-z^2q^2)(1-z^2q^3)}  {\tilde {\cal J}}(zq)= {(1-z^3q^6)\over
(1-zq^2)(1-z^2q^2)(1-z^2q^3)}{\tilde {\cal J}}(zq)\,,
\eeq
whose  solution is (\ref{casdeu}).

The expression (\ref{casdeu}) can be easily expressed as a $q$-series:
\beq
{\tilde {\cal J}}(z)= \sum_{m_1,m_2,m_3\geq0}{(-1)^{m_3}q^{2m_1+2m_2+3m_3(m_3+3)/2}\, 
z^{m_1+2m_2+3m_3}\over (q)_{m_1}(q)_{m_2}(q^3;q^3)_{m_3}}\,.\eeq
Although this is not a manifestly positive $q$-series, the positivity is inherited from the first equality in
(\ref{redt}).
To illustrate this result, we see that the coefficient of $z^5$ is $q^6+2q^7+4q^8+7q^9+\cdots$ and the seven
corresponding
$02$-partitions of weight 9 are
\beq
{{\cal J}}(5,9)=7 \; : \{(50202),\; (41202),\; (32202),\; (31212),\; (31302),\; (23202),\; (22212)\}\,.\eeq

\n {\bf Corollary 22.} The generating function for the number of partitions of $n$ into $m$ parts and subject to the
restriction
\beq \label{diffQa}
\la_i\geq \la_{i+1}\; \qquad {\rm and}\qquad \la_i\geq \la_{i+2}+4\;,\eeq
is
\beq \label{uuut}
P(z)= \sum_{m_1,m_2,m_3\geq0}{(-1)^{m_3}q^{2m_1+2m_2+3m_3(m_3+3)/2+(m_1+2m_2+3m_3)
(m_1+2m_2+3m_3-2)}\,
z^{m_1+2m_2+3m_3}\over (q)_{m_1}(q)_{m_2}(q^3;q^3)_{m_3}}\,.\eeq

\n {\it Proof.} $02$-partitions are transformed into ordinary partitions by adding  the staircase
$(\cdots,9,7,5,3,1,-1)$:
\beq
\la_i= n_i+2(m-i)-1\,. \eeq
It is immediate to verify that the $\la_i$'s satisfy (\ref{diffQa}). Hence, by modifying the power of $q$ within
${\tilde {\cal J}}(z)$ to take into account the addition of the staircase, which amounts to replace $q^nz^m\rw
q^{n+m(m-2)}z^m$, we obtain the generating function $P(z)$ for partitions satisfying (\ref{diffQa}) in the form
(\ref{uuut}).

This expression appears to be a new result. It differs from the generating functions pertaining to the same counting
problem presented in Theorem 5.14 of \cite{Mu3} and Theorem 2.7 of \cite{Mu4}.

There is of course a whole tower of $0p$-partitions, with $p\geq 1$. We have seen that for $p=1$ and $p=2$, the
generating functions for jagged partitions with specified length are solutions of a first-order $q$-equation. This does not
seem to be the case for $p>2$. For instance, the counting of  03-partitions lead to a second-order equation.

\subsection{012-partitions}

We now turn to the analysis of jagged partitions of another type, namely the 012-partitions.
 
\n {\bf Definition 23.} $012$-partitions  are defined as vectors $(n_1,\cdots, n_m) $ of
non-negative entries
$n_j\geq 0$,with
$n_m\geq 2$, subject to the  weak ordering conditions:
\beq \label{weaky}
n_j\geq n_{j+1}-1\;,\qquad  n_j\geq n_{j+2}-2\;,\qquad  n_j\geq n_{j+3}\,.
\eeq

\n {\bf Theorem 24.} The generating function for ${\cal J'}(m,n)$, the number of $012$-partitions of $n$ with length $m$,
is 
\beq \label{casdeuxa}
 {\tilde {\cal J'}}(z)=\sum_{n,m\geq0} {\cal J'}(m,n)z^mq^n= {(-zq^2)_\y (-z^2q^3)_\y\over
 (z^3q^3)_\y}\,.
\eeq

\n {\it Proof. } 
Introduce the following  function-set bijections:
\begin{eqnarray}
& & {\cal J'}(m,n):\, (\cdots 012^+)\;, \quad  {\cal K'}(m,n): \,(\cdots 112^+)\;,\quad 
 {\cal L'}(m,n):\, (\cdots 122^+)\;,\cr
& &{\cal M'}(m,n):\, (\cdots 212^+)\;,\quad{\cal N'}(m,n):  (\cdots
222^+)\,.
\end{eqnarray}
 These are related as follows
\begin{eqnarray}
& & (\cdots 012012^+)=(\cdots 012^+012) + (\cdots 112^+)\cr
& & (\cdots 112112^+)=(\cdots 112^+112) + (\cdots 122^+12)+ (\cdots 122^+) \cr
& & (\cdots 122122^+)=(\cdots 212^+2) + (\cdots 123^+)\cr 
& & (\cdots 212212^+)=(\cdots 122^+12) + (\cdots 222^+)\cr
& & (\cdots 222222^+)=(\cdots 222^+2) + (\cdots 223^+)\,.
\end{eqnarray}
The functional retranscription of these recurrences reads
\begin{eqnarray}
& & {\cal J'}(m,n)= {\cal J'}(m-3,n-3) +{\cal K'}(m,n)\cr
& & {\cal K'}(m,n)= {\cal K'}(m-3,n-4) +{\cal L'}(m-2,n-3)+{\cal L'}(m,n)\cr
& & {\cal L'}(m,n)= {\cal M'}(m-1,n-2) +{\cal J'}(m,n-m)\cr
& & {\cal M'}(m,n)= {\cal L'}(m-2,n-3) +{\cal N'}(m,n)\cr
& & {\cal N'}(m,n)= {\cal N'}(m-1,n-2) +{\cal K'}(m,n-m)\,.\
\end{eqnarray}
leading to 
\begin{eqnarray}
& & {\tilde {\cal J'}}(z)=z^3q^3{\tilde {\cal J'}}(z) +{\tilde {\cal K'}}(z)\cr
& & {\tilde {\cal K'}}(z)= z^3q^4{\tilde {\cal K'}}(z) +z^2q^3{\tilde {\cal L'}}(z)+{\tilde {\cal L'}}(z)\cr
& & {\tilde {\cal L'}}(z)= zq^2{\tilde {\cal M'}}(z) +{\tilde {\cal J'}}(zq)\cr
& & {\tilde {\cal M'}}(z)= z^2q^3{\tilde {\cal L'}}(z) +{\tilde {\cal N'}}(z)\cr
& & {\tilde {\cal N'}}(z)= zq^2{\tilde {\cal N'}}(z) +{\tilde {\cal K'}}(zq)\,.\
\end{eqnarray}
The auxiliary functions  can be eliminated and we get
\beq \label{redta}
{\tilde {\cal J'}}(z)= {(1+z^2q^3)(1-z^4q^8)\over (1-zq^2)(1-z^3q^3)(1-z^3q^4)(1-z^3q^5)}  {\tilde {\cal
J'}}(zq)\,.\eeq for which the solution is (\ref{casdeuxa}).

\n {\bf Corollary 25.} The generating function for the number of partitions of $n$ into $m$ parts and subject to the
restriction
\beq \label{diffQ}
\la_i\geq \la_{i+1}\; \qquad {\rm and} \qquad \la_i\geq \la_{i+3}+3\, ,\eeq
is
\beq \label{uuutr}
P'(z)= \sum_{m_1,m_2,m_3\geq0}{q^{m_1(m_1+3)/2+m_2(m_2+5)/2+3m_3+(m_1+2m_2+3m_3)
(m_1+2m_2+3m_3-3)/2}\, 
z^{m_1+2m_2+3m_3}\over (q)_{m_1}(q)_{m_2}(q)_{m_3}}\,.\eeq

\n {\it Proof. } To generate partitions from 012-partitions of length $m$, we add the staircase $(m-2,\cdots, 1,0,-1)$
of weight $m(m-3)/2$, i.e., $\la_j=n_j-m-j-1$. The weak ordering conditions  (\ref{weaky}) imply the restrictions
(\ref{diffQ}).  Therefore, by replacing
$z^N $ by $ z^N q^{N(N-3)/2}$ in the generating function for
${\cal J'}(m,n)$, whose series expansion reads
 \beq 
{\cal J'}(z)=\sum_{m_1,m_2,m_3\geq0}{q^{m_1(m_1+3)/2+m_2(m_2+5)/2+3m_3}\, 
z^{m_1+2m_2+3m_3}\over (q)_{m_1}(q)_{m_2}(q)_{m_3}}\,,\eeq
we recover the above expression for $P'(z)$.  

This expression $P'(z)$ agrees with the suitable specialization of the generating function displayed in Eq. 1.6 of
\cite{Mu1} (see also Theorem 9.9 of
\cite{Mu2}).

The consideration of 0123-partitions leads to a second-order $q$-equation. We thus expect that counting jagged
partitions with ground state of period
$0123\cdots p$ will always lead to higher-order
$q$-equations for
$p\geq 3$.

\subsection{$001$-partitions }

The final class of jagged partitions to be considered are the $0\cdots 01$-ones (with $p$ zeros), written for short as
$0^p1$-partitions.

\n {\bf Definition 26.} $0^p1$-partitions  are defined as vectors $(n_1,\cdots, n_m) $ of
non-negative entries
$n_j\geq 0$,with
$n_m\geq 1$, subject to the  weak ordering conditions:
\beq 
n_j\geq n_{j+s}-1\;\quad{\rm for}\quad  1\leq s\leq p,\qquad  n_j\geq n_{j+p+1}\,. \
\eeq

For all the cases we have considered (namely, $1\leq p\leq 6$), we always end up with a first-order
$q$-equation, suggesting that this is a generic feature of this class  of jagged partitions. Unfortunately, the
resulting generating functions we obtain do not have a nice form for $p>2$, in particular, we cannot write them as
a multi-sum.
 For this reason, we
limit ourself to displaying the generating function for 001-partitions. Moreover, since its proof is quite similar to the
previous ones, we simply give  the recurrence relation underlying its construction:
\begin{eqnarray}
& & (\cdots 001001^+)= (\cdots 001^+001) +(\cdots 011^+01)+ (\cdots 011011^+)\cr
& & (\cdots 011011^+)=(\cdots 011^+011) + (\cdots 1111^+) \cr
& & (\cdots 111111^+)=(\cdots 111^+1) + (\cdots 112112^+)\,.\
\end{eqnarray}

\n {\bf Theorem 27. }  The generating function for ${\cal J''}(m,n)$, the number of $001$-partitions of $n$ with
length
$m$, is 
\beq \label{casdeux}
 {\tilde {\cal J''}}(z)= {(-z^2q;q^2)_\y \over
 (zq)_\y\, (z^3q;q^3)_\y\, (z^3q^2;q^3)_\y}\,.
\eeq

Again, we can apply this result to the enumeration of partitions subject to some restrictions. In the present case, the
restriction is a superposition of a `difference-1 condition at distance 2' and `difference-3 condition at distance 3'. The
partition $\la=(\la_1,\cdots,\la_m)$  is obtained from the vector $(n_1,\cdots,n_m)$ augmented by the staircase
$(m-1,\cdots,1,0)$.

\n {\bf Corollary 28.} The generating function for the number of partitions of $n$ into $m$ parts and subject to the
restrictions
\beq \label{diffQq}
\la_i\geq \la_{i+1}\; ,\qquad \quad \la_i\geq \la_{i+2}+1  \qquad {\rm and}\qquad \la_i\geq \la_{i+3}+3 \,,\eeq
is
\beq \label{uuutra}
P''(z)= \sum_{m_0,m_1,m_2,m_3\geq0}{q^{\beta}\, 
z^{2m_0+m_1+3m_2+3m_3}\over (q^2;q^2)_{m_0}\, (q)_{m_1}\, (q^3;q^3)_{m_2}\, (q^3;q^3)_{m_3}}\,,\eeq
with
\beq
\beta= m_0^2+m_1+m_2+2m_3+(2m_0+m_1+3m_2+3m_3)(2m_0+m_1+3m_2+3m_3-1)/2\,.
\eeq

\appendix 
 
\section{Proof of identity (\ref{idrty})}
We first obtain an identity for $\zeta_0^2+ 4q\zeta_4^2$. Starting from the expression of the $\zeta_i$'s given in
(\ref{leszeta}), we have
\begin{eqnarray} \label{urr}
\left[{(q^8;q^8)_\y^2\over (q^{16};q^{16})_\y} \right]^2(\zeta_0^2+ 4q\zeta_4^2)\,  &= &
\left[{(q^4;q^4)_\y^{5}\over (q^{2};q^{2})_\y^2(q^{8};q^{8})_\y^2}\right]^2 + q \left[{2(q^8;q^8)_\y^2\over
(q^{4};q^{4})_\y}\right]^2
 \cr
 &= &
\left[\sum_{n=-\y}^\y q^{2n^2}\right]^2 + q \left[\sum_{n=-\y}^\y q^{2n(n+1)}\right]^2
 \cr
 &= &
\sum_{n,m=-\y}^\y [q^{2n^2+2m^2}+ q^{2n^2+2n+2m^2+2m+1}]
 \cr
 &= &
\sum_{n,m=-\y}^\y [q^{(n+m)^2+(n-m)^2}+ q^{(n-m)^2+(n+m+1)^2}]
 \cr
 &= &
\sum_{n,m=-\y}^\y q^{(n-m)^2}[q^{(n+m)^2}+ q^{(n+m+1)^2}]
 \cr
 &= &
\sum_{n,m=-\y}^\y q^{n^2}[q^{(n+2m)^2}+ q^{(n+2m+1)^2}]
 \cr
 &= &
\sum_{n,m=-\y}^\y q^{n^2}q^{(n+m)^2}
 = 
\sum_{n,m=-\y}^\y q^{n^2}q^{m^2}
 \cr
&= &
\left[\sum_{n=-\y}^\y q^{n^2}\right]^2
= 
{(q^2;q^2)_\y^{10}\over (q)_\y^4\, (q^4;q^4)_\y^4}\,.
\end{eqnarray}
The starting trick is to multiply the sum $\zeta_0^2+ 4q\zeta_4^2$ by an appropriate factor to transform it into a
sum of squares. Each square is then rewritten as an infinite sum by means of the Jacobi triple-product
identity. In the fifth line, we replace $n$ by $n-m$. In the sixth one, we see that the second summation is broken into
its even and odd parts; grouping them together leads to the seventh line.  There, we replace $m$ by $m-n$.
At the end, Jacobi triple-product
identity is used once more. The core identity in (\ref{urr}) is of course
\beq   \label{urrg}
\left[\sum_{n=-\y}^\y q^{2n^2}\right]^2 + q \left[\sum_{n=-\y}^\y q^{2n(n+1)}\right]^2= \left[\sum_{n=-\y}^\y
q^{n^2}\right]^2\,,
\eeq
which is actually well-known\footnote{It also has a natural combinatorial interpretation. Write the
$q$-series expansion of (\ref{urrg}) as $\sum_{p=0}^\y a_pq^p$. The right hand side shows that
$a_p$ is the number of ways $p$ can be written as a sum of two squares. The left hand side shows that for
$p$ even, $a_p$ is also the number of ways of writing $p/2$  as a sum of two squares, while for $p$ odd, $a_p$
stands for the number of ways $(p-1)/2$ can be written as a sum of two near squares, i.e., products
$n(n+1)$.}.
With 
\beq
\theta_2(q)= \sum_{n=-\y}^\y q^{(n+1/2)^2}\;, \qquad \theta_3(q)= \sum_{n=-\y}^\y q^{n^2}\;,
\eeq
 (\ref{urrg}) is nothing but
\beq
\theta_2^2(q^2)+\theta_3^2(q^2)= \theta_3^2(q)\;,\eeq
an identity that arises naturally from the parametrization of the
arithmetic-geometric mean iteration using  theta functions - cf. section 2.1 of \cite{Bor} (and 
(\ref{urrg}) is precisely eq. (2.1.8) there). Let us now return to (\ref{urr}) in order to finish the proof of
(\ref{idrty}). Isolating
$\zeta_0^2+ 4q\zeta_4^2$ and multiplying the result by
$\zeta_0\zeta_1^3\zeta_4$ yields:
\beq 
\zeta_0\zeta_1^3\zeta_4( \zeta_0^2+ 4q\zeta_4^2) = {(q^2;q^2)_\y^{14} \, (q^{16};q^{16})_\y^{7}\over (q)_\y^7\,
(q^8;q^8)_\y^{14}}=
\zeta_1^7\,,
\eeq
as desired.

  \vskip0.3cm
\noindent {\bf ACKNOWLEDGMENTS}

We thank the referee for  corrections (in particular, on the range of values of $r$ for which Lemma 13 is correct),
for valuable suggestions and for drawing our attention to the fact that (\ref{urrg}) is a standard theta
function identity. We are also grateful to  O. Warnaar for useful discussions (and in particular, for pointing out the
combinatorial interpretation of (\ref{urrg})) and for his comments on the article.  We also  thank E. Mukhin for
guiding us through the recent literature of
$(k,\ell)$ admissible partitions and Luc B\'egin for his collaboration in
 \cite{BFJM}.
This work is supported by NSERC.

  \vskip0.3cm
\noindent {\bf NOTE ADDED}

A bijection between jagged partitions and overpartitions has also been obtained by K. Mahlburg, {\it The overpartition function modulo small powers of 2}, Discrete Math. {\bf 286} (2004), 263-267. This article  describes an extension of the congruence mod 64 obtained here.


\begin{thebibliography}{99}
\addcontentsline{toc}{section}{References}





\bibitem{Andr}
G.E. Andrews, {\it The theory of
partitions}, Cambridge Univ. Press (1984).



\bibitem{Andrr}
G.E. Andrews, {\it Multiple $q$-series}, Houston J. Math. {\bf 7} (1981) 11-22.

\bibitem{AAR}
G.E. Andrews, R. Askey and R. Roy, {\it Special functions}, Encyclopedia of Mathematics and its applications {\bf 71},
Cambridge Univ. Press (1999).


\bibitem{BFJM}
L. B\'egin, J.-F. Fortin, P. Jacob and P. Mathieu,  {\it Fermionic characters
for graded parafermions}, Nucl. Phys. {\bf B659} (2003) 365-386.

\bibitem{Bor}
J.M. Borwein and P.B. Borwein, {\it Pi and the AGM - a study in analytic number theory and computational complexity},
Wiley, N.Y., (1987).


\bibitem{Mu4}
B. Feigin, M. Jimbo, S. Loktev, T. Miwa, E. Mukhin, {\it  Bosonic formulas for $(k,l)$-admissible partitions},
math.QA/0107054.


\bibitem{Mu3} 
B. Feigin, M. Jimbo, T. Miwa, {\it Vertex operator algebra arising from the minimal series $M(3,p)$ and monomial basis},
math.QA/0012193. 




\bibitem{Mu1}
B. Feigin, M. Jimbo, T. Miwa, E. Mukhin and Y. Takeyama, {\it Particle content of the $(k,3)$-configurations},
math.QA/0212348.


\bibitem{Mu2}
B. Feigin, M. Jimbo, T. Miwa, E. Mukhin and Y. Takeyama, {\it Fermionic formulas for $(k,3)$-admissible configurations},
math.QA/0212347.






\bibitem{FJM}
 J.-F. Fortin, P. Jacob and P. Mathieu,  {\it Generating function 
for $K$-restricted jagged partitions}, math-ph/0305055.



\bibitem{Har}
G.H. Hardy, {\it Ramanujan}, Cambridge Univ. Press (1940).

\bibitem{Rad}
H. Rademacher, {\it Topics in analytic number theory},  Springer Verlag (1973).




\end{thebibliography}
\end{document}